\documentclass[12pt,reqno]{amsart}
\usepackage{amssymb}
\usepackage[all]{xy}

\newtheorem{thm}{Theorem}[section]
\newtheorem{prop}[thm]{Proposition}
\newtheorem{lem}[thm]{Lemma}
\newtheorem{cor}[thm]{Corollary}

\theoremstyle{definition}

\newtheorem*{ack}{Acknowledgment}

\theoremstyle{remark}
\newtheorem*{notat}{Notations and Conventions}

\newtheorem*{remark}{Remark}

\numberwithin{equation}{section}

\newcommand{\cm}{Cohen-Macaulay}
\newcommand{\tr}{\operatorname{tr}}
\newcommand{\ch}{\operatorname{char}}
\newcommand{\height}{\operatorname{height}}
\newcommand{\depth}{\operatorname{depth}}
\newcommand{\res}{\operatorname{res}}
\newcommand{\cores}{\operatorname{cor}}
\newcommand{\Fract}{\operatorname{Fract}}
\newcommand{\rank}{\operatorname{rank}}
\newcommand{\diag}{\operatorname{diag}}
\newcommand{\ann}{\operatorname{ann}}
\newcommand{\Ker}{\operatorname{Ker}}
\newcommand{\GL}{\operatorname{GL}}
\newcommand{\End}{\operatorname{End}}
\newcommand{\Aut}{\operatorname{Aut}}
\newcommand{\I}{\operatorname{\mathbb{I}}}
\newcommand{\Ext}{\operatorname{Ext}}

\newcommand{\tto}{\longrightarrow}

\newcommand{\Mat}{{\operatorname{M}}}
\newcommand{\Mn}{\Mat_n}

\newcommand{\cX}{\mathcal{X}}
\newcommand{\cE}{\mathcal{E}}
\newcommand{\cP}{\mathcal{P}}
\newcommand{\cCM}{\mathcal{C\mspace{-6mu}M}}
\newcommand{\fq}{\mathfrak{q}}
\newcommand{\fQ}{\mathfrak{Q}}
\newcommand{\fa}{\mathfrak{a}}

\newcommand{\bbF}{\mathbb{F}}
\newcommand{\bbN}{\mathbb{N}}
\newcommand{\bbC}{\mathbb{C}}
\newcommand{\bbZ}{\mathbb{Z}}
\newcommand{\bbQ}{\mathbb{Q}}
\newcommand{\Abar}{\overline{A}}
\newcommand{\Ghat}{\widehat{G}}

\begin{document}

%
%

\title[Rings of invariants]{On Cohen-Macaulay rings of invariants}

\author{M. Lorenz}
\address{Department of Mathematics, Temple University,
	Philadelphia, PA 19122-6094}
\email{lorenz@math.temple.edu}

\author{J. Pathak}
\email{pathak@math.temple.edu}
\thanks{Research of both authors supported in part by NSF grant DMS-9988756}

\subjclass{13A50, 16W22, 13C14, 13H10}
\keywords{finite group action, ring of invariants, invariant theory, 
height, depth, Cohen-Macaulay ring, cohomology, Sylow subgroup}

\begin{abstract}
We investigate the transfer of the Cohen-Macaulay property from a  commutative 
ring to a subring of invariants under the action of a finite group. Our point of 
view is ring theoretic and not a priori tailored to a particular type of group action. 
As an illustration, we briefly discuss the special case of multiplicative actions, 
that is, actions on group algebras  $k[\bbZ^n]$  via an action on $\bbZ^n$. 
\end{abstract}

\maketitle

\section*{Introduction}

This article addresses the question to what extent the {\cm} property passes
from a (commutative) ring $R$ to a subring $R^G$ of invariants under the
action of a finite group $G$ on $R$. As is well-known, the {\cm} property is
indeed inherited by $R^G$ whenever the trace map $\tr_G \colon R \to R^G$,
$r\mapsto \sum_{g\in G} g(r)$, is surjective (\cite{hochster}; see also Section
\ref{SS:nonCM} below). In the opposite case, however, the property
usually does not transfer, even in the particular case of linear actions, that is,
$G$-actions on polynomial
algebras $R=k[X_1,\dots,X_n]$ by linear substitutions of the variables. The 
{\cm} problem for linear invariants has been rather thoroughly
explored without, at present, being anywhere near a final solution.

Our focus in this article will not be on linear $G$-actions on polynomial algebras
nor, for the most part, on any other kind of group action on affine algebras
over a field. Rather, in Sections~\ref{S:reltrace} -- \ref{S:sylow}, we work entirely in the setting of commutative
noetherian
 rings. Besides being marginally more general, this approach has resulted in a
number of simplifications of results previously obtained by Kemper \cite{ke1}, \cite{ke2}
in a geometric setting using geometric methods. Nevertheless, the article
owes a great deal to Kemper's insights and originated from a study of his work.

A rough outline of the contents is as follows. Section~\ref{S:reltrace} is devoted to
relative trace maps. We determine the height of their image, an ideal of $R^G$, and use
this result to give a lower bound for the height of annihilators in $R^G$ of 
certain cohomology classes. Section~\ref{S:depth} reviews basic material on {\cm} rings
and local cohomology and describes a pair of spectral sequences constructed by
Ellingsrud and Skjelbred \cite{es}. These are used to derive certain depth estimates. In
Section~\ref{S:cm}, we return to rings of invariants $R^G$ and note some easy facts on the
non-{\cm} locus of $R^G$ and on the special case of Galois actions; it turns out that if
the $G$-action on $R$ is Galois in the sense of Auslander and Goldman \cite{auslander}
then $R^G$ is {\cm} if and only if $R$ is. Section~\ref{S:modules} develops the main technical
tools of this article. We use the aforementioned spectral sequences of
Ellingsrud and Skjelbred to derive a depth formula for modules of invariants
which underlies our subsequent applications. The latter concern the
case where $R$ has characteristic $p$ and focus on the role played by the Sylow
$p$-subgroup of $G$. For the precise statements of these results, we refer the reader
to Section \ref{S:sylow} where they are presented. The final Section~\ref{S:multiplicative} 
initiates the study of
the {\cm} property in the special
case of multiplicative actions. These are defined to be $G$-actions on Laurent 
polynomial algebras
$R= k[X_1^{\pm 1}, \dots,X_n^{\pm 1}]$ stabilizing the lattice of monomials 
$\langle X_1,\dots,X_n\rangle \cong \bbZ^n$; so we may think of $G$ as a 
subgroup of $\GL_n(\bbZ)$.  We show that if $G$ maps onto some non-trivial $p$-group
and has a cyclic Sylow $p$-subgroup, $P$, then $R^G$ is {\cm} if and only if 
$P$ is generated by a bireflection, that is, a matrix $g\in \GL_n(\bbZ)$ so that
$g - 1_{n\times n}$ has rank at most $2$. In this case, $P$ must have order $2$, $3$,
or $4$. A more detailed study of the {\cm} property for multiplicative invariants will
form the subject of the second author's Ph.D.~thesis.

\begin{notat}
Throughout, $G$ will denote a finite group and $R$ will be a commutative
ring on which $G$ acts by ring automorphisms, $r \mapsto g(r)$. The subring
of $G$-invariant elements of $R$ will be denoted by $R^G$ and the skew group
ring of $G$ over $R$ by $RG$. Thus, $RG$ is the free left $R$-module with
basis the elements of $G$, made into a ring by means of the multiplication
rule $rg\cdot r'g'=rg(r')gg'$ for $r,r'\in R$, $g,g'\in G$. The ring $R$
is a module over $RG$ via $rg\cdot r'=rg(r')$.
All modules are understood to be left modules.
\end{notat}


\section{The relative trace map} \label{S:reltrace}

\subsection{} \label{SS:reltraceprelim}
Throughout this section, $H$ denotes a subgroup of $G$.
The \emph{relative trace
map} $\tr_{G/H}: R^H \to R^G$ is defined by
$$
\tr_{G/H}(r) = \sum_{g \in G/H} g(r)  \qquad (r \in R^H)\ .
$$
Here, $g$ runs over any transversal for the cosets $gH$ of $H$ in $G$.
Since $\tr_{G/H}$ is $R^G$-linear, the image of $\tr_{G/H}$ is an ideal of $R^G$ which we
shall denote by
$$
R_H^G\ .
$$

\subsection{Covering primes}
The proof of the following lemma was communicated to us by Don Passman.
The special case where $R$ is an affine algebra over a field is covered by
\cite[Satz 4.7]{ke2}.
As usual, we will write $^gH=gHg^{-1}$ $(g\in G)$ and
$I_G(\fQ) = \{g \in G \mid (g-1)(R) \subset \fQ\}$ denotes the \emph{inertia group} of 
an ideal $\fQ$ of $R$.

\begin{lem}\label{L:primes}
For any prime ideal $\fQ$ of $R$,
$$
\fQ \supseteq R_H^G \iff [I_G(\fQ):I_{^gH}(\fQ)]\in \fQ\quad \text{for all $g\in G$}
$$
\end{lem}

\begin{proof}
The implication $\Leftarrow$ follows from the straightforward formula
\begin{equation*}\label{E:tr}
\tr_{G/H}(r) \equiv \sum_{g\in I_G(H)\backslash G/H} [I_G(\fQ):I_{^gH}(\fQ)]\,g(r) \mod \fQ
\end{equation*}
for all $r\in R^H$. For $\Rightarrow$, assume that $\fQ\supseteq R_H^G$. 
It suffices to show that 
$$
[I_G(\fQ):I_H(\fQ)]\in \fQ \ .
$$
Indeed, $R_H^G = R_{^gH}^G$, since $\tr_{G/H}(r)=\tr_{G/^gH}(g(r))$ holds for all 
$r\in R^H$ and $g\in G$.

To simplify notation, 
put $I=I_G(\fQ)$ and let $P$ denote
a Sylow $p$-subgroup of $I\cap H = I_H(\fQ)$, where  
$p\ge 0$ is
the characteristic of the commutative domain $R/\fQ$. (Here $P=\{1\}$ if $p=0$.) Then our
desired conclusion, $[I:I\cap H] \in \fQ$, is equivalent with 
$$
[I:P] \in \fQ\ .
$$
Furthermore, our assumption $\fQ\supseteq R_H^G$ entails that $\fQ\supseteq R_P^G$, because
$\tr_{G/P}=\tr_{G/H}\circ\tr_{H/P}$. Thus, leaving $H$ for $P$, we may assume that
$H=P$ is a $p$-subgroup of $I$. Let $D=\{g\in G \mid g(\fQ)=\fQ\}$ denote the decomposition
group of $\fQ$; so $I\le D$. 
We claim that
$$
\fQ \supseteq R_P^D\ .
$$
To see this, choose $r \in R$ so that $r \in g(\fQ)$ for all $g \in G
\setminus D$ but $r \notin \fQ$. Then $s = \prod_{g \in D} g(r)$ also belongs to
$\bigcap_{g \in G \setminus D} g(\fQ)$ but not to $\fQ$ 
and, in addition, $s \in R^D$. Now assume that, contrary to our claim, there
exists an element $f \in R^P$ so that $\tr_{D/P}(f) \notin \fQ$. Then
$\tr_{D/P}(s f) = s\tr_{D/P}(f) \in \bigcap_{g \in G \setminus D} g(\fQ) \setminus \fQ$,
and hence $\tr_{G/P}(s f) \notin \fQ$, a contradiction.

By the claim, we may replace $G$ by $D$, thereby reducing 
to the case where $\fQ$ is $G$-stable. (Note that $I$ is unaffected by this replacement.)  
So $G$ acts on $R/\fQ$ with kernel $I$, $P$ is a $p$-subgroup of $I$, 
and $R_P^G\subseteq \fQ$. Thus, $0\equiv \tr_{G/P}(r)\equiv [I:P]\cdot\sum_{g\in G/I} g(r) 
\mod \fQ$
holds for all $r\in R^P$. Our desired conclusion, $[I:P]\in \fQ$, will follow if we
can show that $\sum_{g\in G/I} g(r) \notin \fQ$ holds for some $r\in R^P$. But
$\sum_{g\in G/I} g$ induces a nonzero endomorphism on $R/\fQ$, by linear independence of 
automorphisms of $K=\Fract(R/\fQ)$; so $\sum_{g\in G/I} g(s) \notin \fQ$ holds 
for some $s\in R$. 
Putting $r=\prod_{h\in P} h(s)$, we have $r\in R^P$ and $r \equiv s^{|P|} \mod \fQ$.
Since $|P|$ is $1$ or a power of $p=\ch K$, we obtain 
$\sum_{g\in G/I} g(r) \equiv \sum_{g\in G/I} g(s^{|P|})
\equiv \left(\sum_{g\in G/I} g(s)\right)^{|P|} \notin \fQ$, as required. 
\end{proof}

\subsection{Height formula} \label{SS:height}
For any collection $\cX$ of subgroups of $G$, we define the ideal $R_{\cX}^G$
of $R^G$ by
$$
R_{\cX}^G = \sum_{H\in\cX}R_H^G\ .
$$
Inasmuch as $R_D^G\subseteq R_H^G = R_{^gH}^G$ holds for all $D\le H\le G$ and 
$g\in G$, there is no 
loss in assuming that $\cX$ is closed under $G$-conjugation and under taking subgroups.

Moreover, for any subgroup $H\le G$, we define 
$$
I_R(H) = \sum_{h\in H} (h-1)(R)R\ .
$$
Thus, $I_R(H)$ is an ideal of $R$, and $\fQ\supseteq I_R(H)$ is equivalent with
$H\le I_G(\fQ)$.

\begin{lem} \label{L:height}
Assume that $\bbF_p\subseteq R$, and let $\cX$ be a collection of
subgroups of $G$ that is closed under $G$-conjugation and under taking subgroups.
Then
$$
\height R_{\cX}^G = \inf \{\height I_R(P) \mid 
\text{$P$ is a $p$-subgroup of $G$, $P\notin\cX$}\}\ .
$$
\end{lem}

\begin{proof} One has
$$
\height R_{\cX}^G = \inf_{\fq}\height\fq = \inf_{\fQ}\height \fQ\ ,
$$
where $\fq$ runs over the prime ideals of $R^G$ containing $R_{\cX}^G$
and $\fQ$ runs over the primes of $R$ containing $R_{\cX}^G$. Here, the
first equality is just the definition of height, while the second equality
is a consequence of the standard relations between the primes of $R$ and $R^G$;
see, e.g., \cite[Th\'eor\`eme 2 on p.~42]{bourbaki}.

By Lemma~\ref{L:primes},
$$
\fQ \supseteq R_{\cX}^G \iff p\, \big|\, [I_G(\fQ):I_H(\fQ)]\ \text{for all $H\in\cX$.}
$$
Since $I_H(\fQ)=I_G(\fQ)\cap H$ belongs to $\cX$ for $H\in\cX$, the latter
condition just says that the Sylow $p$-subgroups of $I_G(\fQ)$ do not belong 
to $\cX$ or, equivalently,
some $p$-subgroup $P\le I_G(\fQ)$ does not belong to $\cX$. Therefore,
$$
\fQ \supseteq R_{\cX}^G \iff \fQ \supseteq 
\bigcap_{\text{$P\le G$ a $p$-subgroup, $P\notin\cX$}}I_R(P)\ ,
$$
which implies the asserted height formula.
\end{proof}

\subsection{Annihilators of cohomology classes}\label{SS:ann}
Let $M$ be a module over the skew group ring $RG$. Then, for each 
$r\in R^G$, the map $\rho \colon M \to M$, $m\mapsto rm$, is
$G$-equivariant, and hence $\rho$ induces a map on cohomology $\rho_*: H^*(G,M)
\to H^*(G,M)$. Letting $r$ act on $H^*(G,M)$ via $\rho_*$ we make 
$H^*(G,M)$ into an $R^G$-module.

\begin{lem} \label{L:ann}
The ideal $R_H^G$ of $R^G$ annihilates the kernel of the restriction map
$\res^G_H: H^*(G,M) \to H^*(H,M)$.
\end{lem}

\begin{proof}
The action of $R^G=H^0(G,R)$ on $H^*(G,M)$ can also be interpreted as 
coming from the cup product
$$
H^0(G,R)\times H^*(G,M) \stackrel{\cup}{\tto} H^*(G,R\otimes M)
\stackrel{\cdot}{\tto} H^*(G,M) \ ,
$$
where the map denoted by $\cdot$ comes from the $G$-equivariant 
map $R\otimes M\to M$,
$r\otimes m\mapsto rm$; see, e.g., \cite[Exerc.~1 on p.~114]{brown}.
Furthermore, the relative trace map $\tr_{G/H}: R^H \to R^G$ is identical
with the corestriction map $\cores^G_H: H^0(H,R)\to H^0(G,R)$; 
cf.~\cite[p.~81]{brown}. Thus, the transfer formula for cup products
(\cite[(3.8) on p.~112]{brown}) gives, for $s\in R^H$ and $x\in H^*(G,M)$,
$$
\tr_{G/H}(s)x = \cdot(\tr_{G/H}(s)\cup x) = \cdot(\cores_H^G(s\cup\res^G_H(x)))\ .
$$
Therefore, if $\res^G_H(x)=0$ then $\tr_{G/H}(s)x=0$.
\end{proof}

We summarize the material of this section in the following proposition. For convenience,
we write $\res^G_P(\,.\,) = \,.\,\big|_P$.

\begin{prop} \label{P:height}
Assume that $\bbF_p\subseteq R$, and let $M$ be an $RG$-module. 
Then, for any $x\in H^*(G,M)$,
$$
\height\ann_{R^G}(x) \ge \inf\{\height I_R(P) \mid 
\text{$P$ a $p$-subgroup of $G$, $x\big|_P\neq 0$}\}\ .
$$
\end{prop}

\begin{proof}
Let $\cX$ denote the splitting data of $x$, that is, 
$\cX=\{H\le G \mid x\big|_H=0\}$. By Lemma~\ref{L:ann}, 
$\ann_{R^G}(x) \supseteq R_{\cX}^G$, and by Lemma~\ref{L:height},
$\height R_{\cX}^G = \inf \{\height I_R(P) \mid 
\text{$P$ is a $p$-subgroup of $G$, $x\big|_P\neq 0$}\}$.
The proposition follows.
\end{proof}

\section{Depth} \label{S:depth}

\subsection{}
In this section, $A$ denotes any commutative noetherian ring,
$\fa$ is an ideal of $A$, and $M$ denotes
a finitely generated module over the group ring $A[G]$.

\subsection{Depth and local cohomology}\label{SS:depth}
Let $H^i_{\fa}$ denote the $i$-th local cohomology functor with respect to
$\fa$, that is, the $i$-th right derived functor of the $\fa$-torsion functor
$$
\Gamma_{\fa}(M)= H^0_{\fa}(M)=\{m\in M \mid \text{$m$ is annihilated by some power of
$\fa$}\}\ .
$$
Then
$$
\depth(\fa,M) = \inf\{i \mid H^i_{\fa}(M) \neq 0\}
$$
(where $\inf\emptyset=\infty$); see \cite[Theorem  6.2.7]{brodmann}.

Recall from Section~\ref{SS:ann} (with $A=R^G$) that 
$H^*(G,M)$ is a module over $A$. 
Our hypotheses on $A$ and $M$ entail that $M$ is a noetherian 
$A$-module, and hence so are all $H^q(G,M)$. Therefore, 
$$
\depth(\fa,M) = \inf\{i \mid H^i_{\fa}(M) \neq 0\}
$$
and
$$
\depth(\fa,H^q(G,M)) = \inf\{i \mid H^i_{\fa}(H^q(G,M)) \neq 0\}\ .
$$
All $H^i_{\fa}(M)$ are $A[G]$-modules, via the action of $A[G]$ on $M$.

\subsection{The Ellingsrud-Skjelbred spectral sequences}
The above $A$-modules $H^p_{\fa}(H^q(G,M))$ feature as the $E_2^{pq}$-terms 
of a certain spectral
sequence due to Ellingsrud and Skelbred \cite{es}. In fact, two related
spectral sequences are constructed in \cite{es} in the following manner.
 
The
$\fa$-torsion functor $\Gamma_{\fa}$ and the $G$-fixed point functor 
$(\,.\,)^G=H^0(G,\,.\,)$ clearly
commute: $\Gamma_{\fa}(M^G)=\left(\Gamma_{\fa}(M)\right)^G$. Moreover, if the 
$A[G]$-module $M$ is injective,
then one checks that $\Gamma_{\fa}(M)$ is also injective as $A[G]$-module (as in
\cite[Prop. 2.1.4]{brodmann}) and $M^G$ is injective as $A$-module.
Therefore, $H^i(G,\Gamma_{\fa}(M))=0$ and
$H^i_{\fa}(M^G)=0$ holds for all $i>0$ if $M$ is injective. 
We obtain two Grothendieck spectral sequences converging to
$H_{\fa}^*(G,M):=R^{*}(\Gamma_{\fa}(\,.\,)^G)(M) = R^{*}((\,.\,)^G\Gamma_{\fa})(M)$,
for any $A[G]$-module $M$; e.g., \cite[Theorem 11.38]{rotman}: 
\begin{equation}\label{E:es}
\xymatrix@=1pt@!{
E_2^{p,q}=H_{\fa}^p(H^q(G,M)) \ar@2{->}[dr] \\
&  H_{\fa}^{p+q}(G,M)  \\
\cE_2^{p,q}= H^p(G,H^q_{\fa}(M)) \ar@2{->}[ur] 
}
\end{equation}

\subsection{Depth estimates} 
The depth formulas in Section~\ref{SS:depth} combined with the spectral sequences
\eqref{E:es} yield the following estimates for $\depth(\fa,M^G)$.

\begin{lem}\label{L:estimates}
\begin{enumerate}
\item \textnormal{\textbf{lower bound:}} $\depth(\fa,M^G)\ge \min\{\depth(\fa,M),h_{\fa}+1\}$,
	where $h_{\fa}=\inf_{q>0}\{q+\depth(\fa,H^q(G,M))\}$.
\item \textnormal{\textbf{upper bound:}} Assume that $H^{p_0}_{\fa}(H^{q_0}(G,M))\neq 0$
	for some $p_0\ge 0$, $q_0>0$ with $s=p_0+q_0 < \depth(\fa,M)$. Assume further that
	$H^{s+1-\ell}_{\fa}(H^{\ell}(G,M)) = 0$ holds for $\ell = 1,\dots,q_0-1$ and
	$H^{s-1-\ell}_{\fa}(H^{\ell}(G,M)) = 0$ holds for $\ell>q_0$. Then
	$\depth(\fa,M^G)\le s+1$.
\end{enumerate}
\end{lem}

\begin{proof}
Put $m=\depth(\fa,M)$.
Then $H^q_{\fa}(M)=0$ for $q<m$, and so the $\cE$-sequence in \eqref{E:es} implies that
$H^n_{\fa}(G,M)=0$ for $n<m$. Therefore, the $E$-sequence satisfies
\begin{equation} \label{E:zero}
E^{p,q}_{\infty} = 0 \quad\text{if $p+q < m$.}
\end{equation}
Furthermore, $E_2^{p,0}=H_{\fa}^p(M^G)$; so
$$
\depth(\fa,M^G)=\inf\{p \mid E_2^{p,0}\neq 0\}\ .
$$
Finally,
$$
h_{\fa} = \inf\{p+q \mid q>0, E_2^{p,q}\neq 0\}\ .
$$

To prove (a), assume that $p < \min\{m,h_{\fa}+1\}$. Then $E_{\infty}^{p,0}=0$, 
by \eqref{E:zero}, and
$E_r^{i,j}=0$ for $j>0$, $i+j<p$, $r\ge 2$. Recall that the differential $d_r$ of
$E_r$ has bidegree $(r,1-r)$. Thus, $E_r^{p,0}$ has no nontrivial boundaries
and consists entirely of cycles. This shows that 
$E_2^{p,0}=E_3^{p,0}=\dots=E_{\infty}^{p,0}$, and hence $E_2^{p,0}=0$. Thus, (a)
is proved.

For (b), we check that $E_2^{s+1,0}\neq 0$. Our hypotheses imply that, at position
$(p_0,q_0)$, all incoming differentials $d_r$ $(r\ge 2)$ are $0$ as well as all
outgoing $d_r$ $(r\ge 2, r\neq q_0+1)$. Therefore, $E_{q_0+1}^{p_0,q_0}=E_2^{p_0,q_0}$
and $E_{\infty}^{p_0,q_0}=E_{q_0+2}^{p_0,q_0}=\Ker(d_{q_0+1}^{p_0,q_0})$. The former
implies that $E_{q_0+1}^{p_0,q_0}\neq 0$, by hypothesis in $(p_0,q_0)$, and the latter
shows that $d_{q_0+1}^{p_0,q_0}$ is injective, because $E_{\infty}^{p_0,q_0}=0$
by \eqref{E:zero}. Thus, $d_{q_0+1}^{p_0,q_0}$ embeds $E_{q_0+1}^{p_0,q_0}$ into
$E_{q_0+1}^{s+1,0}$, forcing the latter to be nonzero. Hence, $E_2^{s+1,0}$ is
nonzero as well, as desired.
\end{proof}

\subsection{{\cm} rings}
For any finitely generated $A$-module $V$, one defines $\dim V = \dim(A/\ann_AV)$ and
$$
\height(\fa,V)=\height(\fa+\ann_{A}V/\ann_{A}V)\ ;
$$
so $\dim V = \sup_{\fa}\height(\fa,V)$. Always,
$$
\depth(\fa,V) \le \height(\fa,V)\ ;
$$
see \cite[Exerc.~1.2.22(a)]{bruns}. The $A$-module $V$ is called \emph{\cm} if
equality holds for all ideals $\fa$ of $A$. In order to show that $V$ is {\cm},
it suffices to check that $\depth(\fa,V) \ge \height(\fa,V)$ holds for all
maximal ideals $\fa$ of $A$ with $\fa \supseteq \ann_AV$.

\section{The {\cm} property for invariant rings} \label{S:cm}

\subsection{}
We now return to invariant rings $R^G$. Our main objective is to 
investigate when the {\cm} property passes from $R$ to $R^G$. 
In this section, we record a few elementary observations that are
independent of the local cohomology methods in Section~\ref{S:depth}.

\subsection{The non-{\cm} locus} \label{SS:nonCM}
By definition, the non-{\cm} locus of $R^G$ consists of those prime ideals
$\fq$ of $R^G$ so that the localization $(R^G)_{\fq}$ is not {\cm}. Thus, $R^G$ is
{\cm} if and only if its non-{\cm} locus is empty. Here, we
point out a general bound for the non-{\cm} locus in terms of relative trace maps. 
More detailed results for affine algebras over a field can be found in
\cite[Kapitel 5]{ke2}. Recall the notation $R_{\cX}^G$ from Section~\ref{SS:height}.

\begin{prop} \label{P:nonCM}
Let $\cCM$ denote the set of subgroups $H$ of $G$ so that $R^H$ is {\cm}.
Then, for every prime ideal $\fq$ of $R^G$ so that $\fq \nsupseteq R_{\cCM}^G$, the 
localization $(R^G)_{\fq}$ is {\cm}.
\end{prop}

\begin{proof}
By hypothesis, $\fq \nsupseteq R_H^G$ for some $H \in \cCM$.
Let $R_{\fq}$ denote the localization of $R$ at the multiplicative subset $R^G\setminus\fq$.
Then the $G$-action on $R$ extends to $R_{\fq}$ and $(R_{\fq})^G = (R^G)_{\fq}$;
see \cite[Prop. 23 on p.~34]{bourbaki}. Similarly, $(R_{\fq})^H = (R^H)_{\fq}$; so
$(R_{\fq})^H$ is {\cm}. By choice of $\fq$ the relative trace map
$\tr_{G/H} \colon (R_{\fq})^H \to (R_{\fq})^G$ is onto. Fix an element $c \in (R_{\fq})^H$
so that $\tr_{G/H}(c)=1$ and define $\rho \colon (R_{\fq})^H \to (R_{\fq})^G$ by
$\rho(x)=\tr_{G/H}(cx)$. This map is a ``Reynolds operator", i.e., $\rho$ is
$(R_{\fq})^G$-linear and restricts to the identity on $(R_{\fq})^G$. Since $(R_{\fq})^H$
is integral over $(R_{\fq})^G$, a result of Hochster and Eagon (\cite{hochster} or 
\cite[Theorem 6.4.5]{bruns}) implies that $(R_{\fq})^G$ is {\cm}, which proves the
proposition.
\end{proof}

As an application, we note that if $G$ has subgroups $H_i$ so that each $R^{H_i}$ is
{\cm} and the indices $[G:H_i]$ are coprime in $R^G$ then $R^G$ is {\cm} as well.
Indeed, writing $1=\sum_i [G:H_i]r_i$ with $r_i\in R^G$, we obtain 
$1=\sum_i \tr_{G/H_i}(r_i) \in R_{\cCM}^G$; so the non-{\cm} locus of $R^G$ is empty.

\subsection{Galois actions} \label{SS:galois}
Recall that the $G$-action on $R$ is \emph{Galois}, in the sense of Auslander and Goldman
\cite{auslander}, if every maximal ideal of $R$ has trivial inertia group in $G$.

\begin{prop} \label{P:galois}
If the $G$-action on $R$ is Galois then $R^G$ is {\cm} if and only if $R$ is.
\end{prop}

\begin{proof}
By \cite[Lemma 1.6 and Theorem 1.3]{chase}, the trace map $\tr_{G/1}\colon R \to R^G$ 
is surjective for Galois actions and $R$ is finitely generated
projective as $R^G$-module. Thus, $R$ is faithfully 
flat as $R^G$-module. Moreover,
for any prime $\fQ$ of $R$ and $\fq = \fQ\cap R^G$, the fibre
$R_{\fQ}/{\fq}R_{\fQ}$ has dimension $0$. Therefore, by \cite[2.1.23]{bruns},
$R^G$ is {\cm} if and only if $R$ is. 
\end{proof}

\section{Modules of invariants} \label{S:modules}

\subsection{}
Throughout this section, $R^G$ is assumed noetherian and $\fa$ denotes an
ideal of $R^G$. Moreover, $M$ denotes
an $RG$-module that is finitely generated as $R^G$-module. Our finiteness assumptions hold,
for example, whenever $R$ is an affine algebra over some noetherian subring 
$k\subseteq R^G$ and $M$ is a finitely generated $RG$-module; 
see \cite[Th\'eor\`eme 2 on p.~33]{bourbaki}. 

\subsection{The problem and a sufficient condition} \label{SS:problem}
Assuming $_RM$ to be {\cm}, 
we are interested in the question under what circumstances $_{R^G}M^G$ will be
{\cm} as well.  We remark that $_RM$ is {\cm} if and only if $_{R^G}M$ is; 
see \cite[Proposition 1.17]{ke2}.

For future reference, we note the following simple lemma.

\begin{lem} \label{L:simple}
Assume that $_RM$ is {\cm} and that $\sqrt{\fa}\supseteq\ann_{R^G}M^G$. Then
$\depth(\fa,M) = \height(\fa,M) \ge \height(\fa,M^G)$. 
\end{lem}

\begin{proof}
Note that $\sqrt{\fa} \supseteq \ann_{R^G}M^G
\supseteq \ann_{R^G}M$ entails that 
$\height(\fa,M)\ge \height(\fa,M^G)$. Further, $\height(\fa,M) = \depth(\fa,M)$,
because $_{R^G}M$ is {\cm}. The lemma follows.
\end{proof}

We now give a sufficient condition for $_{R^G}M^G$ 
to be {\cm}.  We note that $\dim{_RM}  = \dim{_{R^G}M}$, by the usual 
relations between the primes of $R$ and of $R^G$. 

\begin{cor} \label{C:sufficient}
Assume that $_RM$ is {\cm}.
If $H^q(G,M)=0$ holds for $0 < q < \dim{_RM} -1$ 
then $_{R^G}M^G$ is {\cm} as well.
\end{cor}

\begin{proof}
Let $\fa$ be an ideal of $R^G$ with $\fa \supseteq \ann_{R^G}M^G$.
Our hypothesis on $H^q(G,M)$ entails that the value of $h_{\fa}$ in Lemma~\ref{L:estimates}
satisfies $h_{\fa} \ge \dim{_RM} - 1$. Also, $\dim_RM = 
\dim{_{R^G}M} \ge \height(\fa,M) \ge \height(\fa,M^G)$,
by Lemma~\ref{L:simple}. Thus, Lemma~\ref{L:estimates}(a) 
gives $\depth(\fa,M^G) \ge \height(\fa,M^G)$, as required.
\end{proof}

\subsection{Depth formula}
In view of Corollary~\ref{C:sufficient}, we may concentrate on the case where
$M$ has non-vanishing positive $G$-cohomology.
The following proposition is a version of
results of Kemper; see \cite[Corollary 1.6]{ke1} and
\cite[Kor.~1.18]{ke2}.

\begin{prop} \label{P:depth}
Assume that $_RM$ is {\cm} and that $\sqrt{\fa}\supseteq\ann_{R^G}M^G$.
Furthermore, assume that, for some $r\ge 0$, 
$H^q(G,M)=0$ holds for $0 < q < r$ but $\fa x=0$
for some $0\neq x \in H^r(G,M)$. Then
$$
\depth(\fa,M^G) = \min\{r+1,\depth(\fa,M)\}\ .
$$
\end{prop}

\begin{remark}
$\height(\fa,M)=\depth(\fa,M)$ holds in the above formula; see
Lemma~\ref{L:simple}.
\end{remark}

\begin{proof}[Proof of Proposition~\ref{P:depth}]
Our hypothesis $\fa x = 0$ for some $0\neq x \in H^r(G,M)$ is equivalent with
$H^0_{\fa}(H^r(G,M))\neq 0$; so $\depth(\fa,H^r(G,M))=0$. 
The asserted equality is trivial for $r=0$, since 
$\depth(\fa,M^G) = \depth(\fa,M) = 0$ holds in this case. Thus we assume that $r>0$.
Then, in the notation of Lemma~\ref{L:estimates}, we have $r=h_{\fa}$,
and part (a) of the lemma gives the inequality $\ge$.

To prove the reverse inequality, note that Lemma~\ref{L:simple} gives
$\depth(\fa,M) \ge \depth(\fa,M^G)$. Therefore, it suffices to show that 
$\depth(\fa,M^G) \le r+1$
if $\depth(\fa,M)>r+1$. For this, we quote Lemma~\ref{L:estimates}(b) 
with $p_0=0$ and $q_0=r$ (so $s=r$). 
\end{proof}

\section{The Sylow subgroup of $G$} \label{S:sylow}

\subsection{}
In this section, $R$ is assumed to be noetherian as $R^G$-module. 
We further assume that 
$\bbF_p\subseteq R$ and we let $P$ denote a Sylow $p$-subgroup of $G$.

\subsection{A necessary condition} \label{SS:necessary}
Put
$$
\mu = \mu(G,R) = \inf\{r>0 \mid H^r(G,R)\neq 0 \}\ .
$$
\begin{prop} \label{P:necessary}
Put
$\cP = \{P'\le P \mid \height I_R(P') \le \mu+1\}$.
If $R$ and $R^G$ are both {\cm} and $\mu < \infty$ then the restriction map
$$
\res^G_{\cP} \colon H^{\mu}(G,R) \to \prod_{P'\in\cP} H^{\mu}(P',R)  
$$
is injective.
\end{prop}

\begin{proof}
Let $0 \neq x \in H^{\mu}(G,R)$ be given and 
put $\fa=\ann_{R^G}(x)$. Then, by Proposition~\ref{P:height},
$$
\height\fa \ge \inf\{\height I_R(P') \mid 
\text{$P'$ a $p$-subgroup of $G$, $x\big|_{P'}\neq 0$}\}\ .
$$
Since $R^G$ is {\cm}, $\height\fa=\depth\fa$. Finally, Proposition~\ref{P:depth} with $M=R$
gives $\depth\fa \le \mu+1$. Thus, there exists a $p$-subgroup
$P'$ of $G$ with $x\big|_{P'} \neq 0$ and $\height I_R(P') \le \mu+1$.
Note that both the condition $x\big|_{P'} \neq 0$ and the value of $\height I_R(P')$
are preserved upon replacing  $P'$ by a conjugate $^g{P'}$ with $g \in G$. Therefore,
we may assume that $P'\in\cP$, which proves the proposition.
\end{proof}

\subsection{Fixed-point-free actions} \label{SS:fpf}
A subgroup $H$ of $G$ is said to act \emph{fixed-point-freely} on $R$ if
$\height I_R(H') \ge \dim R$ holds for all $1 \neq H' \le H$. 

\begin{cor} \label{C:fpf}
Assume that $R$ is {\cm} and that the Sylow $p$-subgroup of $G$ acts
fixed-point-freely on $R$. Then: $R^G$ is {\cm} if and only if $\dim R \le \mu + 1$.
\end{cor}

\begin{proof}
The implication $\Leftarrow$ follows from Corollary~\ref{C:sufficient} with $M=R$. 
For the converse,
let $R^G$ be {\cm} and assume, without loss, that $\mu < \infty$. Then 
Proposition~\ref{P:necessary}
implies that there is a subgroup $1 \neq P' \le P$ with $\height I_R(P') \le \mu + 1$.
On the other hand, by hypothesis on the $G_p$-action, $\height I_R(P') \ge \dim R$; so
$\dim R \le \mu + 1$.
\end{proof}

\subsection{Bireflections} \label{SS:bireflections}
Following \cite{ke2}, we will call an element $g\in G$ a
\emph{bireflection} on $R$ if $\height I_R(\langle g\rangle) \le 2$.

\begin{cor} \label{C:bireflections}
Assume that $R$ and $R^G$ are {\cm}. 
Let $H$ denote the  subgroup of $G$ that 
is generated by all $p'$-elements of $G$ and all bireflections in $P$.
Then $R^G=R_H^G$.
\end{cor}

\begin{proof}
First note that $H$ is a normal subgroup of $G$ and
$G/H$ is a $p$-group.
Thus, if $R^G \neq R_H^G$ or, equivalently, 
$\widehat{H}^0(G/H,R^H)\neq 0$
then also $H^1(G/H,R^H) \neq 0$; see 
\cite[Theorem VI.8.5]{brown}.
In view of the exact sequence
$$
0 \to H^1(G/H,R^H) \tto H^1(G,R) \stackrel{\res^G_H}{\tto} H^1(H,R)
$$ 
(see \cite[35.3]{baba}) we further obtain $H^1(G,R) \neq 0$. Thus, $\mu = 1$ holds in
Proposition~\ref{P:necessary} and every $P'\in\cP$ consists of bireflections.
Therefore, $P'\subseteq H$ and Proposition~\ref{P:necessary} implies that 
$\res^G_H \colon H^1(G,R) \to
H^1(H,R)$ is injective, contradicting the above exact sequence.
Therefore, we must have $R^G=R_H^G$.
\end{proof}

We remark that if $\bbF_p$ is a $G$-module direct
summand of $R$ then the equality $R^G=R_H^G$ forces $G=H$. 

\subsection{The case $|P|=p$} \label{SS:case}
Put
$$
\mu_p(G) = \mu(G,\bbF_p) = \inf\{r>0 \mid H^r(G,\bbF_p)\neq 0 \}\ .
$$
We will determine this number in the case where the order of $G$ is divisible by 
$p$ but not by $p^2$;
in other words, $|P|=p$. As usual $\bbN_G(P)$ and $\bbC_G(P)$ will denote the normalizer
and the centralizer, respectively, of $P$ in $G$. Thus,
$\bbN_G(P)/\bbC_G(P)$ is a subgroup of $\Aut(P)=\Aut(\bbZ/p)\cong \bbF_p^*$, and hence
it is cyclic of order dividing $p-1$. 

\begin{cor} \label{C:case}
Assume that $|P|=p$. Then $\mu_p(G) = 2[\bbN_G(P):\bbC_G(P)] - 1$. 
Moreover, if $\bbF_p$ is a $G$-module direct
summand of $R$ and $R$ and $R^G$ are both {\cm} then
$\height I_R(P) \le 2[\bbN_G(P):\bbC_G(P)]$. 
\end{cor}

\begin{proof}
Put $N=\bbN_G(P)$, $C=\bbC_G(P)$, and $r = 2[N:C]-1$.
In order to prove that $\mu_p(G)=r$, we use the fact that 
$H^*(G,\bbF_p)\cong H^*(P,\bbF_p)^{N/C}$
holds for $*>0$; see \cite[Corollary 3.6.19]{benson}.
If $p=2$ then $N=C$ and so $r=1$. Moreover, $H^*(P,\bbF_p)^{N/C}\cong
H^*(\bbZ/2,\bbF_2)$ equals $\bbF_2$ in
all degrees. This proves the assertion for $p=2$; so we assume $p$ odd from now on.
In this case, $H^*(\bbZ/p,\bbF_p) \cong \bbF_p[v_1,b_2]/(v_1^2,v_1b_2-b_2v_1)$ with
$\deg v_1=1$ and $\deg b_2=2$; see
\cite[Corollary II.4.2]{adem}. Moreover, identifying $\Aut(\bbZ/p)$ with $\bbF_p^*$,
the action of $\Aut(\bbZ/p)$ on $H^*(\bbZ/p,\bbF_p)$ becomes scalar multiplication,
$v_1 \mapsto \ell v_1$, $b_2 \mapsto \ell b_2$, where $\ell \in \bbF_p^*$. Taking
$\ell$ to be a generator for the subgroup of $\bbF_p^*$ corresponding to $N/C$, we
see that 
$$
H^*(P,\bbF_p)^{N/C} \cong \bigoplus_{i\ge 0} \bbF_p b_2^{i[N:C]} \oplus 
\bigoplus_{i>0} \bbF_p v_1b_2^{i[N:C]-1}\ ;
$$
see \cite[p.~104/105]{adem}. The smallest positive degree where $H^*(P,\bbF_p)^{N/C}$
does not vanish is therefore indeed $2([N:C]-1)+1=r$.

Now assume that $\bbF_p$ is a $G$-module direct
summand of $R$ and $R$ and $R^G$ are both {\cm}. 
The former hypothesis implies that $H^r(G,R)\neq 0$ and hence $\mu\le r$.
Moreover, our hypothesis on $|P|$ implies that $\cP\ni P$ holds in 
Proposition~\ref{P:necessary}, because otherwise $\cP$ would consist of the identity
subgroup alone. Therefore, $\height I_R(P)\le
\mu+1\le r+1$, as desired.
\end{proof}

\section{Multiplicative actions} \label{S:multiplicative}

\subsection{} \label{SS:multprelims}
In this section, we focus on a particular type of group action often
called multiplicative actions. These arise from $G$-actions on lattices
$A\cong \bbZ^n$ by extending this action $k$-linearly to the group algebra $R=k[A]
\cong k[X_1^{\pm 1},\dots,X_n^{\pm 1}]$. Here, we assume $k$ to be a field
such that $p=\ch k$ divides the order of $G$; otherwise the invariant
subalgebra $R^G$ would certainly be {\cm} because $R$ is; see Proposition~\ref{P:nonCM}.
There is no loss in assuming $G$ to be faithfully embedded in $\GL(A)\cong \GL_n(\bbZ)$,
and we will do so. The above notations will remain valid throughout this section.

\subsection{}
A subgroup $H\le G$ acts fixed-point-freely on $R$ if and only if no $1\neq h\in H$
has an eigenvalue $1$ on $A$. Furthermore, an element $g\in G$ is a bireflection on $R$
if and only if the endomorphism $g-1\in\End(A)\cong \Mn(\bbZ)$ has rank at most 2. 
Both observations are consequences of the following 

\begin{lem} \label{L:multheight}
For any subgroup $H\le G$, $\height I_R(H) = n - \rank A^H$.
\end{lem}

\begin{proof}
By definition, the ideal $I_R(H)$ of $R$ is generated by the elements $h(a)-a = 
h(a)a^{-1}-1$ for $h\in H$, $a\in A$. Thus, $R/I_R(H)\cong k[A/[H,A]]$,
where we have put $[H,A]=\langle h(a)a^{-1} \mid h\in H, a\in A\rangle \le A$.
Consequently, $\height I_R(H) = \dim R - \dim R/I_R(H) = n - \rank A/[H,A]$. 
Finally, since the group algebra $\bbQ[H]$ is semisimple,
$A\otimes\bbQ = (A^H\otimes\bbQ) \oplus ([H,A]\otimes\bbQ)$; so
$\rank A/[H,A]= \rank A^H$.
\end{proof}

\subsection{}
Since $G$ permutes the $k$-basis $A$ of $R$, the Eckmann-Shapiro Lemma implies
that
$$
H^*(G,R) \cong \bigoplus_{a\in G\backslash A} H^*(G_a,k) \ ,
$$
where $G_a$ denotes the isotropy group of $a$ in $G$. In particular, using the notations
of Sectios~\ref{SS:necessary} and \ref{SS:case}, we have
\begin{equation} \label{E:mu}
\mu = \inf_{a\in A} \mu_p(G_a) \ .
\end{equation}

\subsection{Example: Inversion} \label{SS:inversion}
Let $G=\langle g = -\I_{n\times n}\rangle$ act on $R = k[X_1^{\pm 1},\dots,X_n^{\pm 1}]$
via $g(X_i) = X_i^{-1}$. This action is fixed-point-free. Moreover, assuming
$p=2$, we have $\mu=\mu_2(G)=1$ by \eqref{E:mu}. Therefore, Corollary~\ref{C:fpf}
gives: 
\begin{quote}
$R^G$ is {\cm} if and only if $n\le 2$.
\end{quote}

\subsection{Example: Reflection groups} \label{SS:reflection}
An element $g\in G$ is called a \emph{reflection} on $R$ if 
$\height I_R(\langle g\rangle)\le 1$
or, equivalently, if the endomorphism $g-1\in\End(A)\cong \Mn(\bbZ)$ has rank at most 1; see
Lemma~\ref{L:multheight}.
If $G$ is generated by reflections then $R^G$ is an affine normal semigroup
algebra over $k$; see \cite{lorenz1}. Therefore, $R^G$ is {\cm} in this case, for any
field $k$; see \cite[Theorem 6.3.5]{bruns}. --- This is in contrast with the situation for
finite group actions on polynomial algebras by linear substitutions
of the variables, where (modular) reflection groups need not lead to {\cm}
invariants \cite{nakajima}.

\subsection{Cyclic Sylow subgroups} \label{SS:cyclic}
As before, we let $P$ denote a fixed Sylow $p$-subgroup of $G$. Moreover,
$O^p(G)$ denotes the intersection of all normal subgroups $N$ of $G$ so that
$G/N$ is a $p$-group.

\begin{thm}
Assume that $O^p(G)\neq G$ and that $P$ is cyclic. Then $R^G$ is {\cm}
if and only if $P$ is generated by a bireflection.
In this case, $P$ has order $2$, $3$, or $4$.
\end{thm}

\begin{proof}
Our hypothesis $O^p(G)\neq G$ is equivalent with $\mu_p(G)=1$; so $\mu=1$ holds as
well, by \eqref{E:mu}. Assuming, $R^G$ to be {\cm}, Corollary~\ref{C:bireflections} 
and the subsequent remark imply that $G=H$. Since all $p'$-elements of $G$ belong
to $O^p(G)$, it follows that $G/O^p(G)=P/P\cap O^p(G)$ is generated by the images 
of the bireflections in $P$. Since $P$ is cyclic, it follows that $P$ is generated
by a bireflection. Now, $P$ acts faithfully on the lattice $A/A^P$ of rank at most $2$.
Thus, $P$ is isomorphic to a cyclic $p$-group of $\GL_2(\bbZ)$, and these are
easily seen to have orders $2$, $3$, or $4$.

The converse follows from the more general Lemma below which does not depend on
cyclicity of $P$ or nontriviality of $G/O^p(G)$.
\end{proof}

\begin{lem}
If $\rank A/A^P \le 2$ then $R^G$ is {\cm}.
\end{lem}

\begin{proof}
By Proposition~\ref{P:nonCM}, it suffices to show that $R^P$ is {\cm};
so we may assume that $G=P$ is a $p$-group. Note that $G$ acts faithfully
on $\Abar = A/A^G$. If $G$ acts as a reflection group
on $\Abar$ then it does so on $A$ as well, and hence the invariants $R^G$ will
be {\cm}; see Section~\ref{SS:reflection}. Thus we may assume that
$\Abar$ has rank $2$ and $G$ acts on $\Abar \cong \bbZ^2$ as a non-reflection $p$-group. 
By the well-known 
classification of finite subgroups of $\GL_2(\bbZ)$ (e.g., \cite[2.7]{lorenz2}),
this leaves the cases $G\cong \bbZ/n$ with $n=2,3$ or $4$ to consider.

The cases $n=2$ or $3$ can be dealt with along similar lines. Indeed, for both values
of $n$,
the only indecomposable $G$-lattices, up to isomorphism, are $\bbZ$, $\bbZ[G]$, and
$\bbZ[G]/(\Ghat)$, where $\Ghat = \sum_{g\in G} g$; see 
\cite[Exercise 4 on p.~514/5]{curtis}. Thus, $A \cong \bbZ^m \oplus (\bbZ[G]/(\Ghat))^r
\oplus \bbZ[G]^s$, and $R^G \cong k[B]^G[X_1^{\pm 1}, \dots, X_m^{\pm 1}]$, where
we have put $B = (\bbZ[G]/(\Ghat))^n \oplus \bbZ[G]^r$. Since $R^G$ is {\cm} if and only if 
$k[B]^G$ is, we may assume that $m=0$. Now,  $\Abar \cong (\bbZ[G]/(\Ghat))^{n+r}$; so 
$2 = (r+s)(|G|-1)$. When $n=3$, this leads to either $r=1$, $s=0$ or
$r=0$, $s=1$. In the former case, $\rank A=2$ and so $R^G$ is surely {\cm}, being a normal 
domain of dimension $2$. If $r=0$, $s=1$
then $A$ is a $G$-permutation lattice of rank $3$. Hence,
$R=k[A]$ is a localization of the symmetric algebra $S(A\otimes k)$, and likewise for the 
subalgebras of invariants. Since linear invariants of dimension $\le 3$ are known to be {\cm}
(e.g., \cite{ke2}), $R^G$ is {\cm} in this case as well. For $n=2$, there are three cases
to consider, one of which ($r=2$, $s=0$) leads to an invariant algebra of dimension $2$ which
is clearly {\cm}. Thus, we are left with the possibilities $r=1$, $s=1$ and $r=0$, $s=2$.
Explicitly, after an obvious choice of basis, $G$ acts as one of the following groups on $A$:
\begin{description}
\renewcommand{\arraystretch}{.75}
\item[Case 1] $G_1 = \left<  g_1=
\left(
\begin{array}{c|c}
{\scriptstyle -1} & \\
\hrulefill & \hrulefill \\
& \begin{smallmatrix} 0 & 1\\ 1 & 0 \end{smallmatrix}
\end{array}
\right)
\right>$;
\item[Case 2] $G_2 = \left<  
\left(
\begin{array}{c|c}
\begin{smallmatrix} 0 & 1\\ 1 & 0 \end{smallmatrix}
 & \\
\hrulefill & \hrulefill \\
& \begin{smallmatrix} 0 & 1\\ 1 & 0 \end{smallmatrix}
\end{array}
\right)
\right>$.
\end{description}
For $G_2$, $A\cong \bbZ^4$ is a permutation lattice. Hence,
as above, it suffices to check that the linear invariant algebra $S(V)^G$ for $V=A\otimes k$
is {\cm} which is indeed the case, by \cite{es}, since $\dim V/V^G = 2$. 
For $G_1$, one can proceed as follows: Embed $G_1$ into
$\Gamma = \langle g_1, \diag(-1,1,1)\rangle \cong \bbZ/2 \times \bbZ/2$
and denote the corresponding basis of $A \cong \bbZ^3$ by $\{x,y,z\}$; so $g_1(x)=x^{-1}$,
$g_1(y)=z$, and $g_1(z)=y$. One easily checks
that $R^{\Gamma} = k[\xi,\sigma_1,\sigma_2^{\pm 1}]$, where $\xi=x+x^{-1}$, $\sigma_1 = y+z$,
and $\sigma_2=yz$. Furthermore, $R = k[A] = R^{\Gamma} \oplus xR^{\Gamma} \oplus
yR^{\Gamma} \oplus xyR^{\Gamma}$. With this, the invariant subalgebra $R^{G_1}$ is easily
determined; the result (for $\ch k = 2$) is 
$R^{G_1} = R^{\Gamma} \oplus (xy+x^{-1}z)R^{\Gamma}$ which is indeed {\cm}.
This completes the proof for $G \cong \bbZ/2$ or $\cong \bbZ/3$.

We now sketch the remaining case, $G\cong \bbZ/4$. The action on
$\Abar = A/A^G$ can then be described by
$G\big|_{\Abar} = \left< s = 
\left(\begin{smallmatrix} 0 & -1\\ 1 & 0 \end{smallmatrix}\right)
\right>$; so $\Abar \cong \bbZ[G]/(s^2+1)$.
With this, one calculates $\Ext_G(\Abar,\bbZ) \cong \bbZ/2$. Thus, there is exactly one
(up to isomorphism) non-split extension of $G$-modules $0 \to \bbZ \to U \to \Abar \to 0$.
A suitable module $U$ is $U = \bbZ[G]/(s-1)(s^2+1)$. Furthermore, one calculates
$\Ext_G(U,\bbZ)=0$. Consequently, either $A \cong A^G \oplus \Abar$ or $A \cong
\bbZ^m \oplus U$, and hence either $R^G \cong k[\Abar]^G[A^G]$ which is {\cm} because
$k[\Abar]^G$ has dimension 2, or $R^G \cong k[U]^G[X_1^{\pm 1},\dots,X_r^{\pm 1}]$
which is {\cm} precisely if $k[U]^G$ is. This reduces the problem to the case where $A=U$
which can be handled by direct calculation, taking advantage of the fact that a conjugate
of group $G_1$ is contained in $G$. We leave the details to the reader.
\end{proof}

\begin{ack}
We thank Don Passman for the proof of Lem\-ma~\ref{L:primes} and Gregor Kemper
for his comments on a preliminary version of this article.
\end{ack}


\end{document}